\documentclass{IEEEtran}
\usepackage[utf8]{inputenc}
\usepackage{cite}
\usepackage{tabu}
\usepackage{array}
\newcolumntype{P}[1]{>{\centering\arraybackslash}p{#1}}
\newcolumntype{M}[1]{>{\centering\arraybackslash}m{#1}}
\usepackage{amsmath,amssymb,amsfonts}
\usepackage{algorithm,algpseudocode,float}
\usepackage{lipsum}
\usepackage{graphicx}
\usepackage{subcaption}
\usepackage{multirow}
\usepackage{textcomp}
\usepackage{xcolor}
\usepackage{amsthm}
\usepackage[english]{babel}
\usepackage{float}
\usepackage{algorithm}
\usepackage{algpseudocode}
\def\BibTeX{{\rm B\kern-.05em{\sc i\kern-.025em b}\kern-.08em
    T\kern-.1667em\lower.7ex\hbox{E}\kern-.125emX}}



\begin{document}

\captionsetup[figure]{font={small},labelformat={default},labelsep=period,name={Fig.}}

\title{Enhancing ACPF Analysis: Integrating Newton-Raphson Method with Gradient Descent and Computational Graphs}

\author{Masoud~Barati,~\IEEEmembership{Senior Member,~IEEE}


\thanks{Masoud Barati is with the Department of Electrical and Computer Engineering and Industrial Engineering, Swanson School of Engineering, University of Pittsburgh, PA, USA (email: masoud.barati@pitt.edu). This work was supported by the NSF ECCS Award 1711921.}}

\maketitle

\begin{abstract}
This paper presents a new method for enhancing Alternating Current Power Flow (ACPF) analysis. The method integrates the Newton-Raphson (NR) method with Enhanced-Gradient Descent (GD) and computational graphs. The integration of renewable energy sources in power systems introduces variability and unpredictability, and this method addresses these challenges. It leverages the robustness of NR for accurate approximations and the flexibility of GD for handling variable conditions, all without requiring Jacobian matrix inversion. Furthermore, computational graphs provide a structured and visual framework that simplifies and systematizes the application of these methods. The goal of this fusion is to overcome the limitations of traditional ACPF methods and improve the resilience, adaptability, and efficiency of modern power grid analyses. We validate the effectiveness of our advanced algorithm through comprehensive testing on established IEEE benchmark systems. Our findings demonstrate that our approach not only speeds up the convergence process but also ensures consistent performance across diverse system states, representing a significant advancement in power flow computation.
\end{abstract}


\section{Introduction}
The integration of renewable energy sources such as wind and solar power is revolutionizing the power systems landscape, presenting new challenges that stem from their variable and intermittent nature. The once-predictable flow of electricity is now subject to fluctuations, leading to a power grid that is more dynamic and less predictable than ever before. This transformation calls for advanced computational techniques capable of conducting power flow analysis with greater resilience and adaptability.

Traditional power flow analysis methods, like the Newton-Raphson (NR) technique, have provided reliable solutions for decades. However, the NR method is primarily designed for stable and predictable systems and may struggle with the irregularities introduced by renewable sources. This is primarily because the NR method's success hinges on good initial approximations and conditions that remain close to normal operating ranges. As renewable integration intensifies, these conditions are increasingly difficult to guarantee, leading to potential convergence issues and inaccuracies. 
Power flow analysis, crucial in the power system field, involves solving nonlinear algebraic equations. The Newton-Raphson (NR) method, widely used for its rapid convergence, iteratively updates solutions using the Jacobian's inverse \cite{1}, \cite{2}. However, this method faces challenges in convergence when initial guesses are far from the final solution or the Jacobian matrix becomes problematic during iterations \cite{3}. Various strategies address these issues, such as augmenting system states \cite{4}, exploring polar versus rectangular formulations \cite{5}, refining starting points \cite{6}, \cite{7}, and employing alternate Jacobian approximations \cite{8, 9, 10}. An innovative approach reformulates power flow as an optimization problem, integrating complementarity constraints for PV buses \cite{11, 12, 13, 14, 15}.

The paper \cite{0} presents a cutting-edge algorithm blending projected gradient descent (GD) and Newton-Raphson (NR) methods, uniquely targeting computational challenges in AC power flow (ACPF) problems. This paper presents an expanded version of \cite{0}, incorporating additional sections on Automatic Differentiation and providing detailed examinations of large power system test cases. This novel strategy redefines the ACPF problem as an optimization challenge, allowing for gradient descent steps without requiring Jacobian matrix inversion, a limitation of conventional NR techniques. Projected GD, effective in maintaining constraints, does not inherently circumvent local optima and saddle points, typical of deterministic optimization. The algorithm smartly transitions to NR methods for quicker convergence as it approaches the global optimum.

This complex scenario demands a sophisticated solution adaptable to the dynamic power grid environment. An effective answer is the integration of the NR method with GD and computational graphs. This comprehensive approach combines NR's iterative resolution prowess with GD's adaptive learning strengths—celebrated for its effectiveness in complex, high-dimensional domains like machine learning and AI.

Computational graphs represent a further leap in this integrated method. By mapping the intricate relationships of power system variables as a network of nodes and edges, computational graphs offer a clear visualization of the power flow problem. They simplify the application of both NR and GD by providing a framework for systematic calculations and updates to the system's state, facilitating the management of the non-linearities characteristic of modern power grids.

In this context, computational graphs not only serve as a visual aid but as a foundational tool that transforms the power flow analysis into a more flexible and adaptive process. This allows for the systematic application of GD, which can iteratively adjust the system state by moving against the gradient of the error surface, thus providing a mechanism to overcome the shortcomings of traditional methods.

The convergence of these methods—NR's precision, GD's adaptability, and computational graphs' clarity—creates a powerful toolkit for today's power system analysts. It equips them to tackle the stochastic nature of renewable energy sources and ensures that power flow analysis remains a reliable and insightful process, crucial for the planning and operation of modern, sustainable power systems. 

The paper is structured to first outline the ACPF problem (Section II), describe the computational graph algorithm and NR method (Section III), provide numerical simulations for six test case studies and comparison with existing methods (Section IV), and conclude with insights and findings (Section V). 

\section{Power Flow Equations}
In an electrical network with $n$ nodes, each node, indexed as $k$, possesses a set of electrical properties: a complex voltage including the magnitude voltage $V_k$, and phase angle $\theta_k$, alongside its associated active $P_k$ and reactive $Q_k$ power components. These properties can be collectively represented in vectorial form as $\mathbf{V}=$ $\left(V_1, \ldots, V_n\right), \boldsymbol{\theta}= \left(\theta_1, \ldots, \theta_n\right), \mathbf{P}=\left(P_1, \ldots, P_n\right)$, and $\mathbf{Q}=\left(Q_1, \ldots, Q_n\right)$. The network's admittance matrix is denoted by $\mathbf{Y}$. This allows the encapsulation of the network's power flow equations into a concise notation: $g(\mathbf{V})=\mathbf{P}+j \mathbf{Q}=\operatorname{diag}\left(\mathbf{V} \mathbf{V}^{\dagger} \mathbf{Y}^{\dagger}\right)$, where $(\cdot)^{\dagger}$ signifies the conjugate transpose operation.
In a power network with a total of $n$ nodes, there are a total of $2n$ distinct real equations. These equations are formulated as follows:
\begin{equation}
\begin{aligned}
& P_k(\mathbf{V}, \boldsymbol{\theta})-P_{\text {net}_k}=0, \quad k=\{1, \cdots, n\} \\
& Q_k(\mathbf{V}, \boldsymbol{\theta})-Q_{\text {net}_k}=0, \quad k=\{1, \cdots, n\}.
\end{aligned}
\label{eq: 19}
\end{equation}

The active and reactive power injections at the $k^{\text{th}}$ node are represented by $P_k(\mathbf{V}, \boldsymbol{\theta})$ and $Q_k(\mathbf{V}, \boldsymbol{\theta})$, respectively. These values are calculated based on the voltage magnitudes and angles. The net active and reactive power entering the $k^{\text{th}}$ node are denoted as $P_{\text{net}_k}$ and $Q_{\text{net}_k}$, respectively. These values are obtained by taking the differences between the generated power ($P_{G k}$, $Q_{G k}$) and the power demand ($P_{D k}$, $Q_{D k}$) at the respective node. The formulas for the active and reactive power injections at each node are as follows:
\begin{equation}
\begin{aligned}
& P_k(\mathbf{V}, \boldsymbol{\theta})=V_k \sum_{m \in \mathcal{G}_k} V_m\left(G_{k m} \cos \theta_{k m}+B_{k m} \sin \theta_{k m}\right) \\
& Q_k(\mathbf{V}, \boldsymbol{\theta})=V_k \sum_{m \in \mathcal{G}_k} V_m\left(G_{k m} \sin \theta_{k m}-B_{k m} \cos \theta_{k m}\right)
\end{aligned}
\label{eq: 20}
\end{equation}

Here, $\mathcal{G}_k$ refers to the set of nodes adjacent to the $k^{\text{th}}$ node. The parameters $G_{k m}$ and $B_{k m}$ represent the conductance and susceptance of the transmission line between nodes $k$ and $m$. The term $\theta_{k m}$ is the angular difference between these nodes.

When given a complex load vector $\mathbf{s}$, the ACPF calculation aims to find the magnitude voltage vector $\mathbf{V}$ and phase angle $\boldsymbol{\theta}$ that satisfy (\ref{eq: eq00}).
\begin{equation}
\mathbf{g}(\mathbf{u,V,\boldsymbol{\theta}})=\mathbf{s}.
\label{eq: eq00}
\end{equation}
Where, the $\mathbf{u}$ is the input of the power flow problem; the active power and magnitude voltage of the PV buses; active power and reactive power of the PQ buses. Rather than directly tackling this nonlinear equation, we propose an optimization framework for its resolution. To solve the system of equations $\mathbf{g}(\mathbf{u,V,\boldsymbol{\theta}})=\mathbf{s}$, we use an optimization approach that aims to minimize the error $\epsilon = \mathbf{g}(\mathbf{u, V, \boldsymbol{\theta}}) - \mathbf{s}$. To quantify the error $\epsilon$, we define a least-square loss function (\ref{eq: eqleast}) as follows. 
\begin{equation}
\min_{\mathbf{{V}, {\boldsymbol{\theta}}}} \frac{1}{2}\|\epsilon\|_2^2.
\label{eq: eqleast}
\end{equation}
It can be rewritten as the following equation:
\begin{equation}
\min_{\mathbf{{V},{\theta}}} \frac{1}{2}\|\mathbf{g}(\mathbf{{V,\boldsymbol{\theta}}})-\mathbf{s}\|_2^2=\min _{\mathbf{V,\boldsymbol{\theta}}} \frac{1}{2} \sum_{i=1}^{n}\left(g_i(\mathbf{V,\boldsymbol{\theta}})-s_i\right)^2
\label{eq: eq1}
\end{equation}
The summation in (\ref{eq: eq1}) includes all PQ, PV, and reference buses, unlike conventional ACPF calculation methods, which aim to maintain equality between the number of variables and equations. It is important to clarify that (\ref{eq: eq1}) does not represent an optimal power flow (OPF) problem. Instead, in order to address the issues of infeasibility and solvability in the ACPF analysis, we approach the ACPF by minimizing the least square error. As a result, we still refer to the problem defined in (\ref{eq: eq1}) as an ACPF problem. In scenarios where the ACPF problem is solvable, the ideal outcome for the objective measure is zero. Under these circumstances, there exists an optimal voltage vector denoted as $\mathbf{V^*,\boldsymbol{\theta}}^*$ that satisfies the condition $\mathbf{g}\left(\mathbf{u,V^*,\boldsymbol{\theta}^*}\right)=\mathbf{s}$. Considering the problem's structure as an unconstrained one characterized by a continuously differentiable objective function, the application of gradient descent emerges as an intuitive method for finding a solution. The loss function $\mathcal{L}$ is given in (\ref{eq: eq1.1}) or (\ref{eq: eq1.2}) contains two real and reactive power flow equations.

\begin{equation}
\min_{\mathbf{{V},{\theta}}} \mathcal{L}= \frac{1}{2}\|\mathbf{g_p}(\mathbf{{V,\boldsymbol{\theta}}})-\mathbf{P}\|_2^2+\frac{1}{2}\|\mathbf{g_q}(\mathbf{{V,\boldsymbol{\theta}}})-\mathbf{Q}\|_2^2
\label{eq: eq1.1}
\end{equation}
\begin{equation}
\min _{\mathbf{V,\boldsymbol{\theta}}} \mathcal{L}= \frac{1}{2} \sum_{i=1}^{n}\left(g_{p_i}(\mathbf{V,\boldsymbol{\theta}})-P_i\right)^2+\frac{1}{2} \sum_{i=1}^{n}\left(g_{q_i}(\mathbf{V,\boldsymbol{\theta}})-Q_i\right)^2.
\label{eq: eq1.2}
\end{equation}
The gradient of $\mathcal{L}$ concerning $\mathbf{V,\boldsymbol{\theta}}$ is derived using the chain rule and can be expressed as:
\begin{equation}
\nabla_{\mathbf{V,\theta}} \mathcal{L}=\mathbf{J}^{\top}(\mathbf{g}(\mathbf{V,\boldsymbol{\theta}})-\mathbf{s}),  
\label{eq: eq2}
\end{equation}
where $\mathbf{J}$ represents the Jacobian matrix linked to the real and reactive power flow equations and different types of reference, PV, and PQ nodes. The formula for conventional Gradient Descent (GD) is as follows:
\begin{equation}
\mathbf{V}_{t+1}=\mathbf{V}_t-\eta \nabla_{\mathbf{V}} \mathcal{L}\left(\mathbf{V}_t\right) 
\label{eq: eq3}
\end{equation}
\begin{equation}
\boldsymbol{\theta}_{t+1}=\boldsymbol{\theta}_t-\eta \nabla_{\boldsymbol{\theta}} \mathcal{L}\left(\boldsymbol{\theta}_t\right),
\label{eq: eq4}
\end{equation}
with $t$ representing the iteration stage, and $\eta$ representing the step size or learning rate, which can be either constant or variable. The sets of node indices for the reference bus, PV buses, and PQ buses are denoted as $\mathcal{I}_{\text {ref }}, \mathcal{I}_{\mathrm{PV}}$, and $\mathcal{I}_{\mathrm{PQ}}$, respectively. The equations (\ref{eq: eq3}) and (\ref{eq: eq4}) describe the gradient descent (GD) algorithm, which is a method used to find the minimum of a function, in this case, the loss function \(\mathcal{L}\). In GD updating equations (\ref{eq: eq3}) and (\ref{eq: eq4}), adjustments are made only to the voltage angles $\theta_i$ for $i \notin \mathcal{I}_{\text {ref }}$ and the voltage magnitudes $V_i$ for $i \in \mathcal{I}_{\mathrm{PQ}}$. This selective updating also allows for the setting of specific voltage magnitudes on $\mathrm{PV}$ buses. In the ACPF problem, our goal is to control the state variables \(\mathbf{x}=(\boldsymbol{\theta},\mathbf{V})\) within specified boundaries, \(\mathbf{x}_{min} \le \mathbf{x} \le \mathbf{x}_{max}\). To solve the constrained minimization problem for an additional set \(K\) that represents the constraints on PV and PQ nodes, we employ the Projected Gradient Descent method and update it as follows:
\begin{subequations}
\begin{align}
& \mathbf{V}_{t+1}^\prime = \mathbf{V}_t-\eta \nabla_{\mathbf{V}} \mathcal{L}\left(\mathbf{V}_t\right) \\
& \mathbf{V}_{t+1} \leftarrow \Pi_K\left(\mathbf{V}_{t+1}^{\prime}\right) 
\end{align}
\begin{align}
& \boldsymbol{\theta}_{t+1}^\prime=\boldsymbol{\theta}_t-\eta \nabla_{\boldsymbol{\theta}} \mathcal{L}\left(\boldsymbol{\theta}_t\right) \\
& \boldsymbol{\theta}_{t+1} \leftarrow \Pi_K\left(\boldsymbol{\theta}_{t+1}^{\prime}\right)
\end{align}
\label{eq: eq4.1}
\end{subequations}
where $\Pi_K\left(\mathbf{x}^{\prime}\right):=\arg \min _{\mathbf{x} \in K}\left\|\mathbf{x}-\mathbf{x}^{\prime}\right\|$ is the projection of $\mathbf{x}^{\prime}$ onto the set $K$, when \(K = \{\mathbf{x} \mid \mathbf{x} \in [\mathbf{x}_{\min}, \mathbf{x}_{\max}]\}\).

\subsection{Conditions for Stopping in GD Algorithim} 
The GD algorithm ceases under two specific conditions:
\subsubsection{Global Optimum is Achieved} This happens when the gradient of the loss function with respect to all parameters (\(\mathbf{V}\) and \(\boldsymbol{\theta}\)) is zero. Mathematically, this condition is represented as \(\mathbf{g}(\mathbf{u,V},\boldsymbol{\theta}) - \mathbf{s} = 0\), where \(\mathbf{g}(\mathbf{u,V}, \boldsymbol{\theta})\) is a function of the parameters and \(\mathbf{s}\) is a desired state or target value. When the gradient is zero, it indicates that the parameters are positioned at a minimum of the loss function, hence no further updates are needed because any small perturbation increases the loss.
\subsubsection{Jacobian Becomes Singular} The Jacobian matrix \(\mathbf{J}\) of partial derivatives of \(\mathbf{g}(\mathbf{u,V}, \boldsymbol{\theta})\) becomes singular. This indicates a point where local changes in \(\mathbf{V}\) and \(\boldsymbol{\theta}\) do not affect the output \(\mathbf{g}\), meaning that the gradient descent cannot effectively update the parameters. Additionally, if \(\mathbf{g}(\mathbf{u,V}, \boldsymbol{\theta}) - \mathbf{s}\) falls within the null space of \(\mathbf{J}^\top\), then changes in parameters have no impact on reducing the error term \(\mathbf{g}(\mathbf{u,V}, \boldsymbol{\theta}) - \mathbf{s}\), effectively making further updates futile.

These stopping conditions ensure that the gradient descent algorithm does not update the parameters endlessly and stops when a minimum is reached or when further updates are ineffective at reducing the loss.
The second situation suggests that the iterative values $\mathbf{V}_t$ and $\boldsymbol{\theta}_t$ are stuck in a local minimum or a saddle point. In order to overcome this obstacle, the algorithm needs to deviate from the gradient path (which becomes zero) and follow a different trajectory. However, this change in direction should not be random, but rather strategically chosen. To address the issue of the iterative values $\mathbf{V}_t$ and $\boldsymbol{\theta}_t$ being trapped in a local minimum or at a saddle point, we need a strategy that enables the gradient descent algorithm to escape these suboptimal points. Here, we propose an enhanced algorithm that incorporates a momentum-based approach along with occasional perturbations to the gradients to help escape local minima and saddle points.

\subsection{Enhanced Gradient Descent Algorithm}
This section outlines an enhanced gradient descent algorithm that employs momentum to sustain movement, adaptive learning rates for step size adjustments, and occasional gradient perturbations, all designed to effectively circumvent local minima and saddle points.

Setting the right hyperparameters for an enhanced gradient descent algorithm involving momentum, perturbation, and adaptive learning rates can significantly impact the effectiveness and efficiency of the training process. Here are the applied settings for these hyperparameters based on power system test cases in the simulation part.
\begin{enumerate}
    \item Momentum Coefficient (\(\gamma\)): the main objective is to accelerate gradients vectors in the right directions, thus leading to faster converging. We set this value between $0.9$ and $0.99$. Start with 0.9 and adjust based on the observed oscillations and convergence rate. Higher values mean more weight is given to the past gradients.
    \item Perturbation Probability (\(p\)): the main purpose is to introduce noise into the gradients to help escape local minima or saddle points. We kept it low to avoid too frequent disruptions of the learning process. We started with a small probability such as $0.05$ or $0.1$ and adjust based on whether the algorithm appears to be stuck often.
    \item 3. Perturbation Intensity (\(\sigma\)): the main goal is to determine the magnitude of the noise added when perturbations occur. The intensity should be small relative to the typical size of the gradient updates. We started with a small fraction of the expected average gradient magnitude, such as $0.01$ or $0.02$. We adjusted based on the stability and effectiveness of the escape mechanism.
    \item Adaptive Learning Rate (\(\eta\)): the main destination is to adjust the step size based on optimization conditions to improve convergence. We considered two options: (a) Constant Learning Rate: Simple but might require manual tweaking; (b) Decay Schemes: Reduces learning rate over time (e.g., exponential decay, step decay). We used a hybrid choice which is to start with a higher rate (e.g., $0.01$ or $0.1$) and reduce it by a factor (e.g., decay by $0.1$ every $10$ epochs).
    \item Adaptive Algorithms: We also utilized the Adam, which automatically adjusts the learning rate and integrates well with momentum. For Adam, $0.001$ as an initial learning rate is applied. For adaptive learning rates, we considered implementing a schedule that reduces the learning rate by a certain factor $0.1$ whenever the decrease in loss plateaus for a specified every $10$ number of epochs. This can help in making fine-grained adjustments towards the later stages of training.
\end{enumerate}
This strategic approach to setting and adjusting hyperparameters will help tailor the optimization process to ACPF model, potentially leading to better performance and convergence.

\begin{algorithm}
\caption{Enhanced Gradient Descent}
\begin{algorithmic}[1]
\State \textbf{Initialize parameters}: $\mathbf{V}_0, \boldsymbol{\theta}_0, \eta, \gamma, p, \sigma$
\State Initialize momentum vectors: $\mathbf{m}_{\mathbf{V}} = \mathbf{0}, \mathbf{m}_{\boldsymbol{\theta}} = \mathbf{0}$

\textbf{for}~{$t = 0, 1, 2, \dots$}
    \State Calculate gradients:
    \[
    \mathbf{g}_{\mathbf{V}} = \nabla_{\mathbf{V}} \mathcal{L}(\mathbf{V}_t)
    \]
    \[
    \mathbf{g}_{\boldsymbol{\theta}} = \nabla_{\boldsymbol{\theta}} \mathcal{L}(\boldsymbol{\theta}_t)
    \]

    \State Apply momentum:
    \[
    \mathbf{m}_{\mathbf{V}} = \gamma \mathbf{m}_{\mathbf{V}} + (1 - \gamma) \mathbf{g}_{\mathbf{V}}
    \]
    \[
    \mathbf{m}_{\boldsymbol{\theta}} = \gamma \mathbf{m}_{\boldsymbol{\theta}} + (1 - \gamma) \mathbf{g}_{\boldsymbol{\theta}}
    \]

    \State Introduce gradient perturbation with probability $p$: \textbf{if}~{$\text{random()} < p$}~\textbf{then}
    \[
        \mathbf{m}_{\mathbf{V}}~+= \sigma \cdot \text{normal}(0, 1, \text{size of}~ \mathbf{m}_{\mathbf{V}}) 
    \]
    \[
        \mathbf{m}_{\boldsymbol{\theta}}~+= \sigma \cdot \text{normal}(0, 1, \text{size of} ~\mathbf{m}_{\boldsymbol{\theta}})
    \]
    \State Update parameters:
    \[
    \mathbf{V}_{t+1} = \mathbf{V}_t - \eta \mathbf{m}_{\mathbf{V}}
    \]
    \[
    \boldsymbol{\theta}_{t+1} = \boldsymbol{\theta}_t - \eta \mathbf{m}_{\boldsymbol{\theta}}
    \]
    \State Check for convergence:
    \textbf{if}~{changes in loss or parameters are below threshold or max iterations reached}~\textbf{then}
        \State \textbf{break}
\end{algorithmic}
\end{algorithm}

\section{Computational Graph}
\subsection{What is computational graph?}
A computational graph is a network where each node signifies an arithmetic operation. It is a structural representation used to depict and compute mathematical expressions efficiently. Consider the elementary mathematical formula:
\begin{equation}
p=x+y    
\end{equation}

\begin{figure}[htp]
    \centering
    \includegraphics[width=8.98cm]{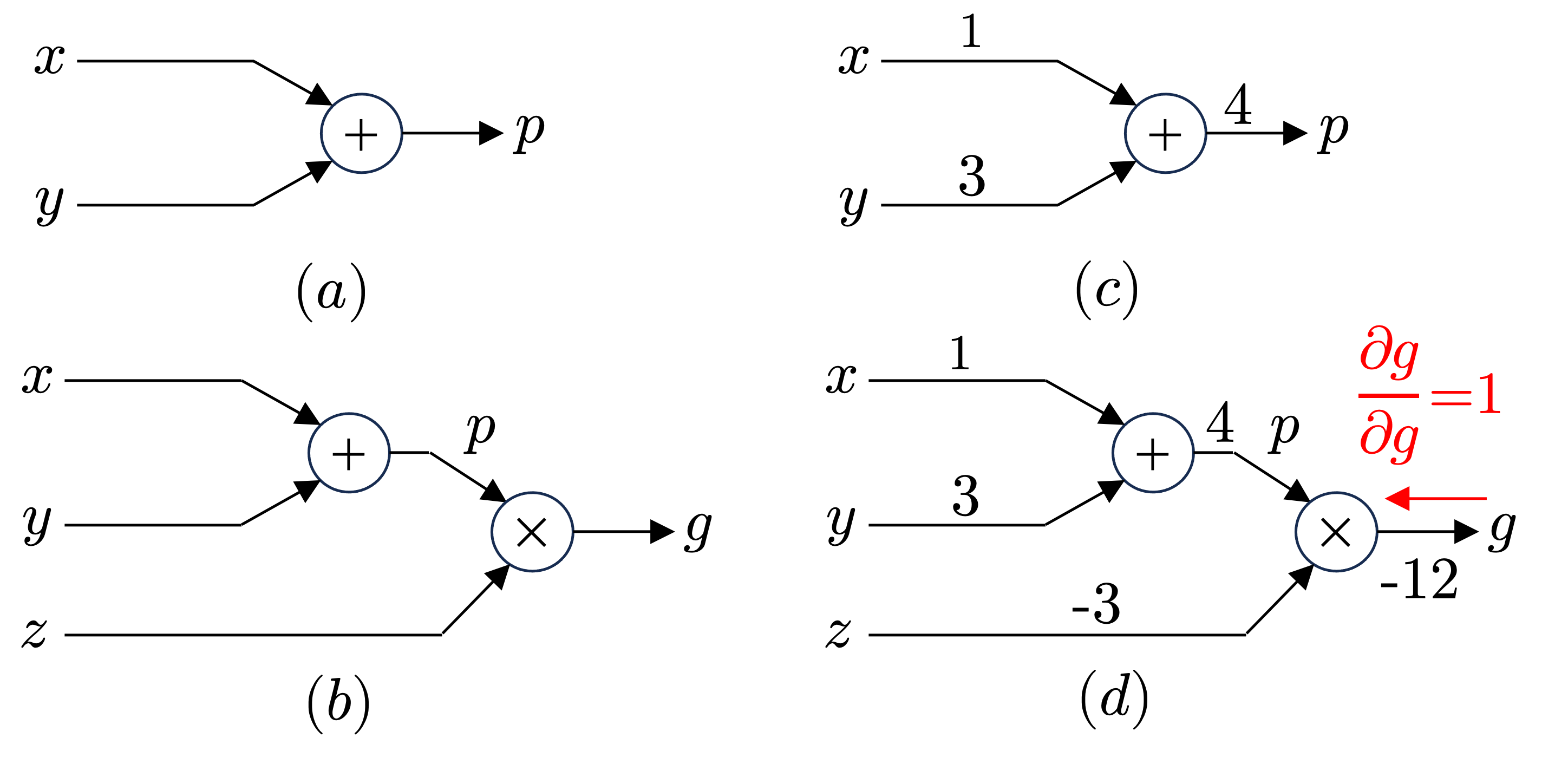}
    \caption{Computational graph for a simple calculation}
    \label{fig1}
\end{figure}

The computational graph for the equation is depicted in Fig. \ref{fig1}(a), where a node marked with a ``+" sign adds inputs $x$ and $y$ to produce output $g$ ". For a complex example, consider the following equation:
\begin{equation}
  g=(x+y) \cdot z.  
\end{equation}

In the computational graph shown in Fig. \ref{fig1}(b), the initial node combines $x$ and $y$ by addition, which then merges with $z$ in a multiplication node to produce the output $g$.

\subsection{Gradient calculation on computational graphs}
In stochastic gradient descent with Newton-Raphson, computational graphs guide each iteration's solution refinement. The forward pass processes inputs sequentially through the graph, moving from initial to terminal nodes, akin to a journey from origin to destination. This involves specific input values progressing through layers and functions at each iteration. For illustrative purposes, let's assign specific values to the inputs as follows: $x=1,~y=3,~z=-3$, as illustrated in Fig. \ref{fig1}(c) and Fig. \ref{fig1}(d). With these values assigned, executing a forward pass allows the computation of intermediary and final output values at each node. In Fig. \ref{fig1}(c), commencing with the values $x=1$ and $y=3$, we calculate the intermediate output $p=4$.
Subsequently in Fig. \ref{fig1}(d), employing $p=4$ and $z=-3$, we determine the final output $g=-12$. The computation progresses linearly from inputs to outputs, with gradient calculation determining each input's impact on the final output, crucial for refining solutions via gradient descent optimization. Consider the necessity to determine the following gradient values: 
\begin{equation}
\frac{\partial x}{\partial f}, \frac{\partial y}{\partial f}, \frac{\partial z}{\partial f}    
\end{equation}

Initiating the backward pass, we calculate the rate of change of the final output in relation to itself, which, by definition, yields a value of one.
\begin{equation}
\frac{\partial g}{\partial g}=1.
\end{equation}

Visualizing our computational graph post this step is shown in Fig.~\ref{fig1}(d). 

Proceeding, we reverse through the multiplication operation. Here, we need to compute the gradient at the nodes $p$ and $z$m where $g$ is the product of $p$ and $z$ ($p=x+y$ and $g=p \cdot z$). 

Using a computational graph in ACPF calculation offers several benefits:
\begin{enumerate}
    \item Efficient Backpropagation: It allows for automatic differentiation, making the calculation of gradients for optimization algorithms (like gradient descent) more efficient and accurate.
    \item Improved Performance: It enables optimization of computational resources and parallel processing, speeding up calculations.
    \item Easier Debugging and Visualization: It helps in visualizing and understanding complex operations, which aids in debugging and improving ACPF models in different use cases and test cases.
\end{enumerate}

Not using a computational graph can lead to:
\begin{enumerate}
    \item Manual Gradient Calculation: This can be error-prone and computationally intensive, especially for complex models.
    \item Reduced Efficiency: Without the optimized execution paths that computational graphs provide, computations may be slower and less efficient.
    \item Difficulty in Scaling: Manual implementations without computational graphs can be challenging to scale for large power grids.
\end{enumerate}

\section{Automatic Differentiation}
This section introduces the concept of automatic differentiation (AD) and classifies methods for computing derivatives in computer programs. The focus is on comparing manual, numerical, symbolic, and automatic differentiation, highlighting the advantages of AD for efficiently calculating derivatives without relying on derivative expressions or suffering from the inaccuracies of numerical approximation.

\subsection{Automatic differentiation differs from alternative differentiation approaches}
This section addresses common misconceptions about automatic differentiation (AD), clarifying that it is distinct from both numerical and symbolic differentiation. AD is described as providing numerical values of derivatives using symbolic rules, thereby combining aspects of both methods without their typical disadvantages. The authors emphasize that AD operates by tracking derivative values during code execution, which allows for precise derivative calculations with minimal overhead.

\subsubsection{Automatic differentiation is not a numerical differentiation} The authors explain that unlike numerical differentiation, which approximates derivatives using finite differences and suffers from accuracy issues due to round-off and truncation errors, AD avoids these problems. Numerical differentiation is straightforward but becomes inaccurate and computationally expensive as the number of dimensions increases, making it unsuitable for applications that need gradients of functions with many variables, such as ACPF problem. Numerical differentiation estimates derivatives by calculating the finite differences at selected sample points of the function. For a function \( P_m: \mathbb{R}^n \rightarrow \mathbb{R} \), the gradient \(\nabla P_m=\left(\frac{\partial P_m}{\partial \theta_2}, \ldots, \frac{\partial P_m}{\partial \theta_n} \mid \frac{\partial P_m}{\partial V_2}, \ldots, \frac{\partial P_m}{\partial V_n}\right)\) can be approximated as:
\begin{equation}
\frac{\partial P_m(\mathbf{x})}{\partial x_i} \approx \frac{P_m\left(\mathbf{x}+h \mathbf{e}_i\right)-P_m(\mathbf{x})}{h},
\end{equation}
where, \(\mathbf{x}=(\boldsymbol{\theta},\mathbf{V})\), \(\mathbf{e}_i\) represents the \(i\)-th unit vector and \(h > 0\) is a small increment. This method is straightforward but requires \(O(n)\) function evaluations for an \(n\)-dimensional gradient, and the choice of \(h\) demands careful attention to avoid inaccuracies.

Numerical methods for derivatives are prone to instability and errors, particularly from the discretization and the inherent limitations of computing precision. As \(h\) decreases, the truncation error diminishes but the round-off error intensifies and can dominate the results.

To reduce such errors, the center difference method is often used:
\begin{equation}
\frac{\partial P_m(\mathbf{x})}{\partial x_i} = \frac{P_m\left(\mathbf{x}+h \mathbf{e}_i\right) - P_m\left(\mathbf{x}-h \mathbf{e}_i\right)}{2h} + O(h^2),
\end{equation}
which balances out first-order errors, improving accuracy by moving the error to the second order in \(h\). Although this method is equally costly as the forward difference in one-dimensional cases, requiring only two function evaluations, it becomes more demanding in ACPF calculation as the number of function dimensions increases, particularly when calculating a Jacobian matrix for functions from \(\mathbb{R}^n\) to \(\mathbb{R}^m\), necessitating \(2mn\) evaluations.

\subsubsection{Automatic differentiation is not a symbolic differentiation} This subsection differentiates AD from symbolic differentiation, which manipulates mathematical expressions to derive formulas for derivatives. Symbolic differentiation can lead to expression swell, where derivatives become unwieldy large expressions, making them difficult to compute and understand. AD, on the other hand, calculates derivatives using actual numerical values during program execution, providing the benefits of the precision of symbolic differentiation without its complexity and inefficiency. Symbolic differentiation automates the process of deriving derivatives from mathematical expressions. This involves transforming expressions using established differentiation rules, such as:
\begin{equation}
\frac{d}{dx_i}(p_{km}(\mathbf{x})+p_{kn}(\mathbf{x})) \rightarrow \frac{d}{dx_i} p_{km}(\mathbf{x})+\frac{d}{dx_i} p_{kn}(\mathbf{x})
\end{equation}
\begin{equation}
\frac{d}{dx_i}(f(\mathbf{x}) g(\mathbf{x})) \rightarrow\left(\frac{d}{dx_i} f(\mathbf{x})\right) g(\mathbf{x})+f(\mathbf{x})\left(\frac{d}{dx_i} g(\mathbf{x})\right)
\end{equation}

In modern computing, tools like Mathematica, Maxima, Maple, and Theano implement this by treating formulas as data structures. The Julia uses structs format. In optimization tasks, symbolic derivatives are crucial for understanding problem structures and can directly provide solutions for extrema, such as finding points where $\frac{d}{dx_i} f(\mathbf{x})=0$, thus bypassing the need for further derivative computations. However, symbolic derivatives tend to become substantially larger than their original expressions, leading to what is known as expression swell—where the size of derivative expressions grows exponentially, complicating their computation. ``Expression swell'' typically refers to the amplification or increase in complexity or size of a mathematical expression in ACPF equations and its derivatives. In this context, it implies that the derivatives of \(P_1\) with respect to \(V_1\) may lead to more intricate or expanded expressions, and the Table \ref{tab:swell} aims to illustrate how these expressions are simplified or reduced in complexity with fewer sinusoidal functions. 

\begin{table}[ht]
\centering
\caption{The derivatives of \(P_1\) with respect to \(V_1\) showcase how changes in \(V_1\) affect the overall expression's swell. In the context of the expression \(P_1^{\prime}(V_1)\), \(\phi_1\) and \(\phi_2\) represent the phase angles associated with the complex admittances \(G_{12} + jB_{12}\) and \(G_{13} + jB_{13}\), respectively.}
\label{tab:swell}
\begin{tabular}{c}
\hline
Original form of $\frac{\partial P_1}{\partial V_1}$                    \\ \hline
$\frac{\partial P_1}{\partial V_1} = 2 G_{11} V_1+V_2\left(G_{12} \cos \theta_{12}+B_{12} \sin \theta_{12}\right)$ \\ ~~~~~~~ $+V_3\left(G_{13} \cos \theta_{13}+B_{13} \sin \theta_{13}\right)$ \\ \hline \hline
Simplified form of $\frac{\partial P_1}{\partial V_1}$ \\ \hline
$\frac{\partial P_1}{\partial V_1} = 2 G_{11} V_1+\sqrt{V_2^2\left(G_{12}^2+B_{12}^2\right)+V_3^2\left(G_{13}^2+B_{13}^2\right)}$ \\ ~~~~~~~~ 
$\cdot \left[\cos \left(\theta_{12}+\phi_1\right)+\cos \left(\theta_{13}+\phi_2\right)\right]$ \\ \hline
\end{tabular}
\end{table}

Considering the function $h(\mathbf{x})=f(\mathbf{x}) g(\mathbf{x})$; both $h(\mathbf{x})$ and its derivative share common elements like $f(\mathbf{x})$ and $g(\mathbf{x})$. Deriving $f(\mathbf{x})$ symbolically and inserting this derivative separately can lead to redundant calculations for any overlapping computations in $f(\mathbf{x})$ and its derivative, resulting in unnecessarily large and complex expressions.
To address this, AD simplifies the process by storing only the necessary intermediate values and interleaving differentiation with simplification steps. This approach, particularly in its forward accumulation mode, optimizes the differentiation process by maintaining the computational efficiency and managing the scale of derivative calculations.

\subsection{Automatic differentiation and its main modes}
This part dives into the technical foundations of AD, particularly detailing its two primary modes: forward mode and reverse mode accumulation. The forward mode computes derivatives alongside the evaluation of the function, which is straightforward but can be computationally expensive for functions with many inputs. Reverse mode, used extensively in ACPF computational graph via backpropagation, calculates derivatives more efficiently by working backwards from the function outputs, making it preferable for functions with a large number of inputs and a smaller number of outputs. 

\begin{figure}[htp]
    \centering
    \includegraphics[width=4.98cm]{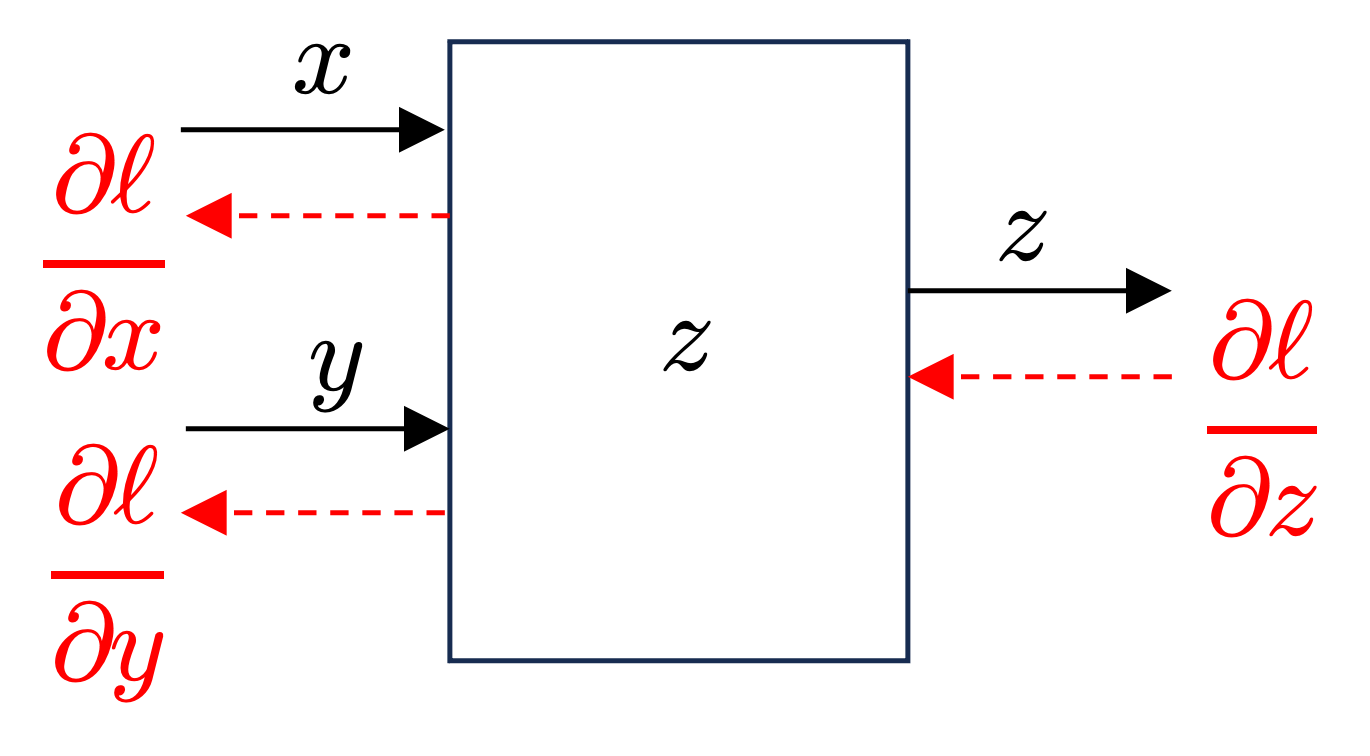}
    \caption{Computational graph illustrating backward differentiation paths for the function $z=$ $z(x,y)$.}
    \label{fig1}
\end{figure}

\begin{figure*}
    \centering
    \includegraphics[width=0.9\linewidth]{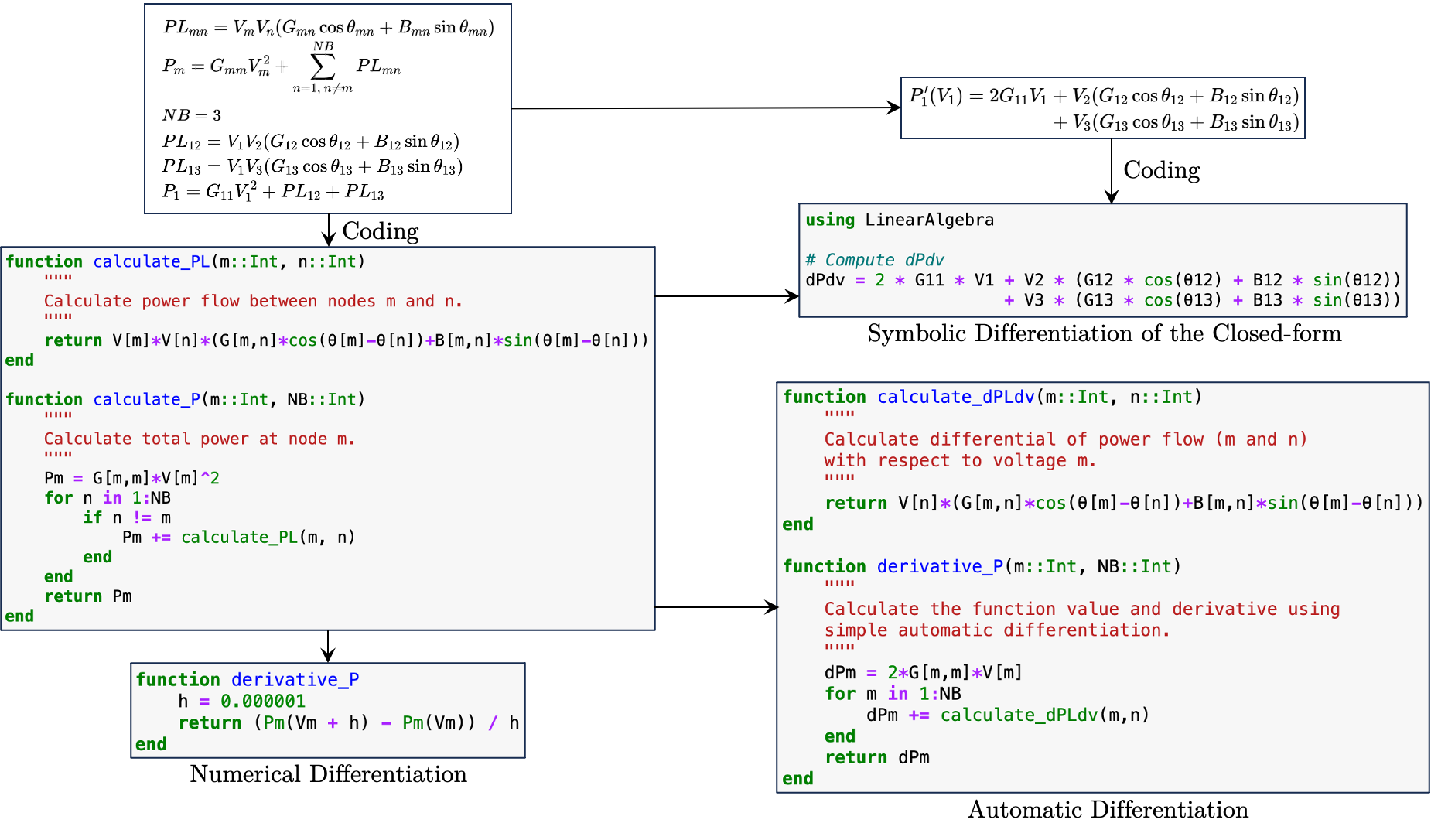}
    \caption{Differentiation techniques for mathematical functions and coding vary in complexity and accuracy. Symbolic differentiation demands exact expressions and can produce complex results, whereas numerical differentiation, simpler but less precise, often encounters errors from data approximations. Automatic differentiation achieves the precision of symbolic methods with less computational overhead and supports dynamic computational features.}
    \label{fig:AD}
\end{figure*}

Figure \ref{fig:AD} illustrates different approaches to obtaining derivatives from mathematical formulas and computational code, as shown by using a simplified logistic map illustration (top left). Symbolic differentiation (middle right) delivers precise outcomes, but requires expressions in a closed format, and struggles with increasing complexity in resulting expressions. Numerical differentiation (bottom right) faces issues with precision due to rounding and truncation errors. However, automatic differentiation (bottom left) achieves accuracy comparable to symbolic methods while maintaining efficiency and allowing for control structures.

Our approach will employ automatic differentiation to determine $\partial z / \partial x$ and $\partial z / \partial y$. Initially, we consider a single node defined by the equation $z=z(x,y)$, which is a component of a broader graph culminating in a scalar value, denoted as $\ell$. Presuming we have successfully computed $\partial \ell / \partial z$, the task then is to find $\partial \ell / \partial x$ and $\partial \ell / \partial y$. To determine the derivative with respect to $x$, the following formula is utilized:
$$
\frac{\partial \ell}{\partial x}=\frac{\partial \ell}{\partial z} \cdot \frac{\partial z}{\partial x}
$$

Similarly, for the derivative with respect to $y$, the equation is:
$$
\frac{\partial \ell}{\partial y}=\frac{\partial \ell}{\partial z} \cdot \frac{\partial z}{\partial y}
$$

It is important to clarify that the expression $\ell(z(x, y))$ may seem slightly confusing. It simply signifies that $\ell$ is a function of $z$, which in turn is a function of $x$ and $y$. In a more complex graph, $\ell$ could depend on numerous other variables.

In computational graphs, understanding the derivative of the final output $\ell$ with respect to a node's output allows reverse calculation of $\ell$'s derivative relative to the node's inputs. This principle enables backpropagation from the final output to initial inputs, a core concept in ACPF calculation. 

\subsection{Gradient Calculation in ACPF}
By reformulating the power flow equation as indicated in equation (\ref{eq: 20}), a more streamlined nested model emerges. This model presents an ideal mathematical structure for the implementation of automatic partial derivative calculations.
\begin{equation}
\begin{aligned}
& P_k(\mathbf{PL})= PL_{kk} + \sum_{m \in \mathcal{G}_k} PL_{km} \\
& Q_k(\mathbf{QL})= QL_{kk} + \sum_{m \in \mathcal{G}_k} QL_{km}
\end{aligned}
\label{eq: 21}
\end{equation}
\begin{equation}
\begin{aligned}
& PL_{km}(\mathbf{c},\mathbf{s})= G_{k m} c_{k m}+B_{k m} s_{k m} \\
& QL_{km}(\mathbf{c},\mathbf{s})= G_{k m} s_{k m}-B_{k m} c_{k m}
\end{aligned}
\label{eq: 22}
\end{equation}
\begin{equation}
\begin{aligned}
& c_{km}(\mathbf{V},\boldsymbol{\theta})= V_k V_m \cos \theta_{k m} \\
& s_{km}(\mathbf{V},\boldsymbol{\theta})= V_k V_m \sin \theta_{k m}
\end{aligned}
\label{eq: 23}
\end{equation}

This revised approach aligns with the computational graph framework and facilitates the automatic computation of gradients. It necessitates the use of a chain rule for the calculation of the Jacobian matrix, integral to gradient determination. This structure is consistent with modern methods of automatic gradient computation, streamlining the process. The power flow equation can be reformulated into a series of nested functions, each dependent on nested variables.
\begin{equation}
\begin{aligned}
& P= PL(c(V,\theta),s(V,\theta))  \\
& Q= QL(c(V,\theta),s(V,\theta)) 
\end{aligned}
\label{eq: 24}
\end{equation}
The gradient can be determined by applying the chain rule in the following manner:

\begin{equation}
\begin{aligned}
& \frac{\partial P}{\partial \theta} = 
\frac{\partial P}{\partial PL} \cdot \frac{\partial PL}{\partial c} \cdot \frac{\partial c}{\partial \theta} + \frac{\partial P}{\partial PL} \cdot \frac{\partial PL}{\partial s} \cdot \frac{\partial s}{\partial \theta}  \\
& \frac{\partial Q}{\partial \theta} = 
\frac{\partial Q}{\partial QL} \cdot \frac{\partial QL}{\partial c} \cdot \frac{\partial c}{\partial \theta} + \frac{\partial Q}{\partial QL} \cdot \frac{\partial QL}{\partial s} \cdot \frac{\partial s}{\partial \theta}  \\
& \frac{\partial P}{\partial V} = 
\frac{\partial P}{\partial PL} \cdot \frac{\partial PL}{\partial c} \cdot \frac{\partial c}{\partial V} + \frac{\partial P}{\partial PL} \cdot \frac{\partial PL}{\partial s} \cdot \frac{\partial s}{\partial V}  \\
& \frac{\partial Q}{\partial V} = 
\frac{\partial Q}{\partial QL} \cdot \frac{\partial QL}{\partial c} \cdot \frac{\partial c}{\partial V} + \frac{\partial Q}{\partial QL} \cdot \frac{\partial QL}{\partial s} \cdot \frac{\partial s}{\partial V}  \\
\end{aligned}
\label{eq: 25}
\end{equation}

Transforming equations (\ref{eq: 20}) and (\ref{eq: 26}) through a first-order Taylor series approximation centered at the point $\left(V^0, \boldsymbol{\theta}^0\right)$ results in the following equations.
\begin{equation}
\begin{aligned}
& P_k\left(\mathbf{V}^0+\mathbf{\Delta V}, \boldsymbol{\theta}^0+\mathbf{\Delta \theta}\right) \\
& =P_k\left(\mathbf{V}^0, \boldsymbol{\theta}^0\right)+\left[\frac{\partial P_k(\mathbf{V}, \boldsymbol{\theta})}{\partial \boldsymbol{\theta}} \mid \frac{\partial P_k(\mathbf{V}, \boldsymbol{\theta})}{\partial \mathbf{V}}\right]_{\left(\mathbf{V}^0, \mathbf\theta^0\right)}\left[\begin{array}{l}
\mathbf{\Delta \theta} \\
\mathbf{\Delta V}
\end{array}\right] \\
& Q_k\left(\mathbf{V}^0+\mathbf\Delta V, \boldsymbol{\theta}^0+\mathbf{\Delta \theta}\right) \\
& =Q_k\left(\mathbf{V}^0, \boldsymbol{\theta}^0\right)+\left[\frac{\partial Q_k(\mathbf{V}, \boldsymbol{\theta})}{\partial \boldsymbol{\theta}} \mid \frac{\partial Q_k(\mathbf{V}, \boldsymbol{\theta})}{\partial \mathbf{V}}\right]_{\left(\mathbf{V}^0, \boldsymbol{\theta}^0\right)}\left[\begin{array}{c}
\mathbf{\Delta \theta} \\
\mathbf{\Delta V}
\end{array}\right] \\
&
\end{aligned}
\label{eq: 26}
\end{equation}

Replacing equations (\ref{eq: 26}) into (\ref{eq: 19}) yields, 
\begin{equation}
\begin{aligned}
& P_{netk}-P_k\left(\mathbf{V}^0, \boldsymbol{\theta}^0\right)=\Delta P_k \\
& =\left[\frac{\partial P_k(\mathbf{V}, \theta)}{\partial \boldsymbol{\theta}} \mid \frac{\partial P_k(\mathbf{V}, \theta)}{\partial \mathbf{V}}\right]_{\left(\mathbf{V}^0, \boldsymbol{\theta}^0\right)}\left[\begin{array}{l}
\mathbf{\Delta \theta} \\
\mathbf{\Delta V}
\end{array}\right] \\
& Q_{n e t k}-Q_k\left(\mathbf{V}^0, \boldsymbol{\theta}^0\right)=\Delta Q_k \\
& =\left[\frac{\partial Q_k(\mathbf{\mathbf{V}, \boldsymbol{\theta})}}{\partial \boldsymbol{\theta}} \mid \frac{\partial Q_k(\mathbf{V}, \boldsymbol{\theta})}{\partial \mathbf{V}}\right]_{\left(\mathbf{V}^0,  \boldsymbol{\theta}^0\right)}\left[\begin{array}{c}
\mathbf{\Delta \theta} \\
\mathbf{\Delta V}
\end{array}\right] \\
&
\end{aligned}
\label{eq: 27}
\end{equation}
From the complete set of linearized power balance equations, we choose a specific subset defined by Equation (\ref{eq: 27}). This subset considers the unique characteristics of each node in the network. It consists of $n$ active power equations, each representing a different node and including the reference node to ensure a complete representation of the Jacobian matrix $\mathbf{J}$. It is important to note that in our proposed model using GD, we do not need to calculate the inverse of the $\mathbf{J}$ matrix like in the traditional Newton-Raphson method. Therefore, we must have a complete representation of the Jacobian matrix in (\ref{eq: eq2}) in order to update the GD and enhanced-GD algorithm 1. The complete form of the Jacobian matrix is established as follows:

\begin{equation}
\left[\begin{array}{l}
\mathbf{\Delta P}_{ref} \\
\mathbf{\Delta P}_{P V} \\
\mathbf{\Delta P}_{P Q} \\
\mathbf{\Delta Q}_{ref} \\
\mathbf{\Delta Q}_{P Q}
\end{array}\right]=\mathbf{J}\left[\begin{array}{l}
0 \\
\mathbf{\Delta \boldsymbol{\theta}}_{P V} \\
\mathbf{\Delta \boldsymbol{\theta}}_{P Q} \\
0 \\
\mathbf{\Delta V}_{P Q}
\end{array}\right]
\end{equation}
where, 
\begin{equation}
\mathbf{J}=\left[\begin{array}{lllll}
0 & \mathbf{H}_{ref, P V} & \mathbf{H}_{ref, P Q} & 0 & \mathbf{N}_{ref, P Q} \\
0 & \mathbf{H}_{P V, P V} & \mathbf{H}_{P V, P Q} & 0  & \mathbf{N}_{P V, P Q} \\
0 & \mathbf{H}_{P Q, P V} & \mathbf{H}_{P Q, P Q} & 0  & \mathbf{N}_{P Q, P Q} \\
0 & \mathbf{J}_{ref, P V} & \mathbf{J}_{ref, P Q} & 0  & \mathbf{L}_{ref, P Q} \\
0 & \mathbf{J}_{P Q, P V} & \mathbf{J}_{P Q, P Q} & 0  & \mathbf{L}_{P Q, P Q}
\end{array}\right]
\end{equation}

\section{Simulation Results}
This section presents the numerical results obtained from implementing the techniques described in the previous sections. The experiments were conducted on Google Cloud Services using the NVIDIA V100, utilizing both CPU and GPU instances for enhanced-gradient descent calculations. We implemented the ACPF model using Julia programming. We examined the efficiency of the enhanced-GD ACPF algorithm for tow set of case studies: Case 1: small scale studies; and Case 2: large scale studies. 

\subsection{Case 1: Small Scale Study}

\begin{table}[ht!]
\centering
\caption{Summary of Initial Solution Violations and Their Magnitudes in Various Test Systems}
\label{tab:my-table}
\begin{tabular}{|c|c|c|}
\hline
\begin{tabular}[c]{@{}c@{}}Case \end{tabular} &
  \begin{tabular}[c]{@{}c@{}}Number of violations \\ (Init. Sol.)\end{tabular} &
  \begin{tabular}[c]{@{}c@{}}Magnitude of the largest \\ violation (p.u.)\end{tabular} \\
 & $V^m$ & $V_{PQ_i}^m$ \\ \hline
IEEE14                            & 7     & 0.0292       \\ \hline
IEEE30                            & 18    & 0.0826       \\ \hline
NEGL39                            & 28    & 0.0603       \\ \hline
IEEE57                            & 36    & 0.0634       \\ \hline
PRCT89                            & 56    & 0.0475       \\ \hline
IEEE118                           & 3     & 0.0070       \\ \hline
\end{tabular}
\end{table}

We used the IEEE 14, 30, 39, 57, 118, and PRCT 98 bus test systems in our tests. The first step involved running an ACPF, and then addressing any potential PQ bus voltage magnitude violations in the load flow solution. To provoke these violations, we raised the lower limits of dependent variables, creating a narrowly feasible region. The results, shown in Tables I and II, demonstrate the effective reduction of violations through orthogonal projections onto feasible regions. The number of violations ranged from 3 in the 118 bus system to 56 in the 89 bus system.

Considering the limitations in control variables (32 for the largest system), it is not practical to address all violations simultaneously. Therefore, only the most significant violations, related to bus magnitude limits, are included in the active constraint set. This approach keeps the system size manageable compared to the original problem and ensures that the necessary reactive power adjustments are minimal yet sufficient to rectify the violations. As observed in Table II's last two columns, this typically results in at least one PQ bus voltage reaching its limit.

\begin{table}[ht]
\centering
\caption{Extended Analysis of Computational Time and Voltage Limits at buses. This table presents the total computational time and the number of buses that meet their voltage limits ($V_i=V_i^{\text {lim}}$) in the final solution for various test systems. The ``CG'' refers to the Proposed Computational Graph method.}
\label{tab:my-table}
\begin{tabular}{|c|c|cc|}
\hline
Case & \begin{tabular}[c]{@{}c@{}}Total (CPU) time\\ (sec) \\ CG / traditional NR\end{tabular} & \multicolumn{2}{c|}{No. of Buses with $V_i=V_i^{\text {lim}}$} \\ \cline{3-4} 
        &                 & \multicolumn{1}{c|}{$\bar{V}_{PQ_i}^m$} & $V_{PQ_i}^M$ \\ \hline
IEEE14  & 0.0025 / 0.0011 & \multicolumn{1}{c|}{2}                  & 0            \\ \hline
IEEE30  & 0.0032 / 0.0028 & \multicolumn{1}{c|}{1}                  & 0            \\ \hline
NEGL39  & 0.4231 / 0.6624 & \multicolumn{1}{c|}{2}                  & 2            \\ \hline
IEEE57  & 0.6443 / 0.8216 & \multicolumn{1}{c|}{1}                  & 1            \\ \hline
PRCT89  & 1.1046 / 1.6421 & \multicolumn{1}{c|}{1}                  & 1            \\ \hline
IEEE118 & 1.3592 / 1.7380 & \multicolumn{1}{c|}{3}                  & 0            \\ \hline
\end{tabular}
\end{table}

The CPU time for computing orthogonal projections is strongly influenced by the number of constraints that are activated or violated. Correcting deviations in the lower voltage limits of the PQ buses may cause the upper voltage limits to be exceeded. The procedure is considered complete only when all such discrepancies are resolved or when it is determined that finding a feasible solution is not possible. In Table II, the second column clearly demonstrates the superiority of the computational graph method over the traditional Newton-Raphson method in large-scale applications. For example, in the case of a 118-bus test system, a significant reduction of 21.79\% in CPU time was observed.

\subsection{Case 2: Large Scale Study}

The terms ``iter'', ``$\frac{d}{dx}$ comp. time'', ``Linear algebra comp. time'', and ``total'' in the table refer to the number of iterations and computational time taken by the algorithm for derivative calculations in automatic differentiation and numerical differentiation, linear algebra computations, and the total processing time, respectively. These metrics are common in optimization contexts to evaluate the performance of algorithms, especially when comparing computational efficiency on different hardware platforms like GPUs. 

Based on the Table IV, which compares two versions of the Enhanced Gradient Descent algorithm—one utilizing auto-differentiation on GPUs and the other using numerical differentiation on GPUs—it is evident that there are significant performance disparities between the two setups:

\begin{enumerate}
     \item Enhanced-GD with Auto-diff (GPU): This method generally shows considerably faster performance in terms of iteration counts and computation times for derivative calculations, linear operations, and overall totals. This enhancement is presumably due to the utilization of GPU acceleration, which is known for its ability to handle parallel computations effectively. For instance, in cases like 500 nodes, the GPU method completes its operations significantly quicker with fewer iterations, demonstrating the computational efficiency of GPUs in handling large-scale calculations.
     \item Enhanced-GD with Numerical-diff (GPU): Although more iterations and longer computation times are required compared to the GPU-based method, it remains a robust choice for environments where GPU resources might not be available. In some extensive cases like 30000 nods and 78484 nods, the GPU version shows a substantial increase in total computation time, which could be a critical factor when dealing with very large datasets or complex computations.
     \item Performance Comparison: In smaller cases, the performance difference is noticeable but not overly dramatic, suggesting that for less demanding tasks, either approach could be suitable depending on the availability of hardware resources. In larger cases, particularly noticeable in 78484 nodes, the Auto-diff approach outperforms the Numerical-diff method dramatically, underscoring the advantage of Auto-diff approach for handling computationally intensive tasks. The Auto-diff method not only completes iterations faster but also manages derivative and linear computations much more efficiently.
\end{enumerate}

\begin{table*}[ht!]
\centering
\caption{The simulation results of enhanced-GD ACPF for automatic differentiation and numerical differentiation}
\label{tab:my-table}
\begin{tabular}{|c|cccc|ccll|}
\hline
\multirow{2}{*}{\textbf{Case}} &
  \multicolumn{4}{c|}{\textbf{Enhanced-GD with Auto-diff}} &
  \multicolumn{4}{c|}{\textbf{Enhanced-GD with Numerical-diff}} \\ \cline{2-9} 
 &
  \multicolumn{1}{c|}{\# iter} &
  \multicolumn{1}{c|}{$\frac{d}{dx}$ comp. time} &
  \multicolumn{1}{c|}{linear alg. comp. time} &
  total &
  \multicolumn{1}{c|}{\# iter} &
  \multicolumn{1}{c|}{$\frac{d}{dx}$ comp. time} &
  \multicolumn{1}{l|}{linear alg. comp.time} &
  total \\ \hline
98 &
  \multicolumn{1}{c|}{32} &
  \multicolumn{1}{c|}{0.02} &
  \multicolumn{1}{c|}{0.04} &
  0.18 &
  \multicolumn{1}{c|}{41} &
  \multicolumn{1}{c|}{0.01} &
  \multicolumn{1}{l|}{0.08} &
  0.08 \\ \hline
179 &
  \multicolumn{1}{c|}{40} &
  \multicolumn{1}{c|}{0.04} &
  \multicolumn{1}{c|}{0.08} &
  0.27 &
  \multicolumn{1}{c|}{64} &
  \multicolumn{1}{c|}{0.01} &
  \multicolumn{1}{l|}{0.14} &
  0.15 \\ \hline
500 &
  \multicolumn{1}{c|}{46} &
  \multicolumn{1}{c|}{0.06} &
  \multicolumn{1}{c|}{0.13} &
  0.46 &
  \multicolumn{1}{c|}{56} &
  \multicolumn{1}{c|}{0.03} &
  \multicolumn{1}{l|}{0.29} &
  0.32 \\ \hline
793 &
  \multicolumn{1}{c|}{44} &
  \multicolumn{1}{c|}{0.03} &
  \multicolumn{1}{c|}{0.07} &
  0.32 &
  \multicolumn{1}{c|}{52} &
  \multicolumn{1}{c|}{0.04} &
  \multicolumn{1}{l|}{0.35} &
  0.39 \\ \hline
1354 &
  \multicolumn{1}{c|}{63} &
  \multicolumn{1}{c|}{0.07} &
  \multicolumn{1}{c|}{0.26} &
  0.71 &
  \multicolumn{1}{c|}{69} &
  \multicolumn{1}{c|}{0.10} &
  \multicolumn{1}{l|}{0.97} &
  1.06 \\ \hline
2312 &
  \multicolumn{1}{c|}{54} &
  \multicolumn{1}{c|}{0.05} &
  \multicolumn{1}{c|}{0.28} &
  0.81 &
  \multicolumn{1}{c|}{64} &
  \multicolumn{1}{c|}{0.14} &
  \multicolumn{1}{l|}{1.57} &
  1.71 \\ \hline
2000 &
  \multicolumn{1}{c|}{48} &
  \multicolumn{1}{c|}{0.05} &
  \multicolumn{1}{c|}{0.14} &
  0.61 &
  \multicolumn{1}{c|}{60} &
  \multicolumn{1}{c|}{0.18} &
  \multicolumn{1}{l|}{1.77} &
  1.94 \\ \hline
3022 &
  \multicolumn{1}{c|}{61} &
  \multicolumn{1}{c|}{0.06} &
  \multicolumn{1}{c|}{0.38} &
  1.04 &
  \multicolumn{1}{c|}{81} &
  \multicolumn{1}{c|}{0.25} &
  \multicolumn{1}{l|}{2.70} &
  2.95 \\ \hline
2742 &
  \multicolumn{1}{c|}{264} &
  \multicolumn{1}{c|}{0.32} &
  \multicolumn{1}{c|}{0.67} &
  3.33 &
  \multicolumn{1}{c|}{148} &
  \multicolumn{1}{c|}{0.69} &
  \multicolumn{1}{l|}{7.88} &
  8.58 \\ \hline
2869 &
  \multicolumn{1}{c|}{70} &
  \multicolumn{1}{c|}{0.07} &
  \multicolumn{1}{c|}{0.16} &
  0.84 &
  \multicolumn{1}{c|}{83} &
  \multicolumn{1}{c|}{0.29} &
  \multicolumn{1}{l|}{3.11} &
  3.40 \\ \hline
3970 &
  \multicolumn{1}{c|}{142} &
  \multicolumn{1}{c|}{0.15} &
  \multicolumn{1}{c|}{2.21} &
  4.26 &
  \multicolumn{1}{c|}{92} &
  \multicolumn{1}{c|}{0.50} &
  \multicolumn{1}{l|}{6.82} &
  7.34 \\ \hline
4020 &
  \multicolumn{1}{c|}{69} &
  \multicolumn{1}{c|}{0.07} &
  \multicolumn{1}{c|}{0.48} &
  1.55 &
  \multicolumn{1}{c|}{84} &
  \multicolumn{1}{c|}{0.49} &
  \multicolumn{1}{l|}{9.55} &
  10.04 \\ \hline
4917 &
  \multicolumn{1}{c|}{68} &
  \multicolumn{1}{c|}{0.07} &
  \multicolumn{1}{c|}{0.47} &
  1.32 &
  \multicolumn{1}{c|}{91} &
  \multicolumn{1}{c|}{0.48} &
  \multicolumn{1}{l|}{5.57} &
  6.05 \\ \hline
4601 &
  \multicolumn{1}{c|}{94} &
  \multicolumn{1}{c|}{0.09} &
  \multicolumn{1}{c|}{0.79} &
  2.07 &
  \multicolumn{1}{c|}{102} &
  \multicolumn{1}{c|}{0.62} &
  \multicolumn{1}{l|}{8.18} &
  8.79 \\ \hline
4837 &
  \multicolumn{1}{c|}{64} &
  \multicolumn{1}{c|}{0.07} &
  \multicolumn{1}{c|}{0.34} &
  1.31 &
  \multicolumn{1}{c|}{83} &
  \multicolumn{1}{c|}{0.55} &
  \multicolumn{1}{l|}{6.69} &
  7.24 \\ \hline
4619 &
  \multicolumn{1}{c|}{61} &
  \multicolumn{1}{c|}{0.08} &
  \multicolumn{1}{c|}{0.22} &
  1.38 &
  \multicolumn{1}{c|}{70} &
  \multicolumn{1}{c|}{0.48} &
  \multicolumn{1}{l|}{8.31} &
  8.78 \\ \hline
5658 &
  \multicolumn{1}{c|}{56} &
  \multicolumn{1}{c|}{0.06} &
  \multicolumn{1}{c|}{0.36} &
  1.34 &
  \multicolumn{1}{c|}{71} &
  \multicolumn{1}{c|}{0.57} &
  \multicolumn{1}{l|}{7.64} &
  8.22 \\ \hline
7336 &
  \multicolumn{1}{c|}{63} &
  \multicolumn{1}{c|}{0.07} &
  \multicolumn{1}{c|}{0.44} &
  1.68 &
  \multicolumn{1}{c|}{70} &
  \multicolumn{1}{c|}{0.69} &
  \multicolumn{1}{l|}{9.91} &
  10.60 \\ \hline
10000 &
  \multicolumn{1}{c|}{75} &
  \multicolumn{1}{c|}{0.09} &
  \multicolumn{1}{c|}{0.55} &
  4.19 &
  \multicolumn{1}{c|}{118} &
  \multicolumn{1}{c|}{1.34} &
  \multicolumn{1}{l|}{17.51} &
  18.86 \\ \hline
8387 &
  \multicolumn{1}{c|}{126} &
  \multicolumn{1}{c|}{0.15} &
  \multicolumn{1}{c|}{0.54} &
  3.12 &
  \multicolumn{1}{c|}{111} &
  \multicolumn{1}{c|}{1.37} &
  \multicolumn{1}{l|}{16.81} &
  18.19 \\ \hline
9591 &
  \multicolumn{1}{c|}{77} &
  \multicolumn{1}{c|}{0.10} &
  \multicolumn{1}{c|}{1.09} &
  3.28 &
  \multicolumn{1}{c|}{98} &
  \multicolumn{1}{c|}{1.34} &
  \multicolumn{1}{l|}{29.74} &
  31.08 \\ \hline
9241 &
  \multicolumn{1}{c|}{575} &
  \multicolumn{1}{c|}{1.80} &
  \multicolumn{1}{c|}{4.79} &
  15.24 &
  \multicolumn{1}{c|}{102} &
  \multicolumn{1}{c|}{1.30} &
  \multicolumn{1}{l|}{18.14} &
  19.44 \\ \hline
1019 &
  \multicolumn{1}{c|}{67} &
  \multicolumn{1}{c|}{0.11} &
  \multicolumn{1}{c|}{0.91} &
  2.93 &
  \multicolumn{1}{c|}{83} &
  \multicolumn{1}{c|}{1.39} &
  \multicolumn{1}{l|}{22.93} &
  24.32 \\ \hline
10480 &
  \multicolumn{1}{c|}{78} &
  \multicolumn{1}{c|}{0.10} &
  \multicolumn{1}{c|}{0.70} &
  3.11 &
  \multicolumn{1}{c|}{98} &
  \multicolumn{1}{c|}{1.43} &
  \multicolumn{1}{l|}{31.68} &
  33.11 \\ \hline
13659 &
  \multicolumn{1}{c|}{88} &
  \multicolumn{1}{c|}{0.13} &
  \multicolumn{1}{c|}{0.91} &
  3.45 &
  \multicolumn{1}{c|}{92} &
  \multicolumn{1}{c|}{1.60} &
  \multicolumn{1}{l|}{22.48} &
  24.08 \\ \hline
20758 &
  \multicolumn{1}{c|}{229} &
  \multicolumn{1}{c|}{0.50} &
  \multicolumn{1}{c|}{6.01} &
  14.44 &
  \multicolumn{1}{c|}{71} &
  \multicolumn{1}{c|}{2.28} &
  \multicolumn{1}{l|}{41.78} &
  44.06 \\ \hline
19402 &
  \multicolumn{1}{c|}{83} &
  \multicolumn{1}{c|}{0.16} &
  \multicolumn{1}{c|}{0.86} &
  5.07 &
  \multicolumn{1}{c|}{101} &
  \multicolumn{1}{c|}{3.16} &
  \multicolumn{1}{l|}{85.12} &
  88.28 \\ \hline
24464 &
  \multicolumn{1}{c|}{75} &
  \multicolumn{1}{c|}{0.16} &
  \multicolumn{1}{c|}{1.03} &
  5.06 &
  \multicolumn{1}{c|}{88} &
  \multicolumn{1}{c|}{2.95} &
  \multicolumn{1}{l|}{57.79} &
  60.75 \\ \hline
30000 &
  \multicolumn{1}{c|}{205} &
  \multicolumn{1}{c|}{0.62} &
  \multicolumn{1}{c|}{4.19} &
  10.98 &
  \multicolumn{1}{c|}{214} &
  \multicolumn{1}{c|}{8.65} &
  \multicolumn{1}{l|}{130.23} &
  138.88 \\ \hline
78484 &
  \multicolumn{1}{c|}{120} &
  \multicolumn{1}{c|}{0.70} &
  \multicolumn{1}{c|}{5.01} &
  23.43 &
  \multicolumn{1}{c|}{145} &
  \multicolumn{1}{c|}{30.10} &
  \multicolumn{1}{l|}{484.30} &
  514.00 \\ \hline
\end{tabular}
\end{table*}
 
\begin{figure}
    \centering
    \includegraphics[width=1\linewidth]{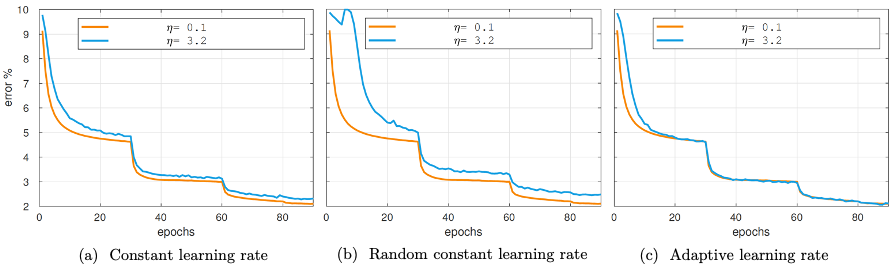}
    \caption{MSE error of enhanced-GD for different learning rate settings.}
    \label{fig:error}
\end{figure}

Figure \ref{fig:error} illustrates the impact of different learning rate strategies on the error rate during training, comparing constant, random constant, and adaptive learning rates over epochs. The constant learning rate shows gradual improvement, with the higher rate (\(\eta = 3.2\)) experiencing some instability. Random constant learning rates introduce more fluctuations, especially at higher rates, suggesting reduced stability due to random changes. The adaptive learning rate strategy outperforms the others by quickly and smoothly reducing error, indicating that it effectively harnesses the benefits of higher rates without the associated drawbacks, thus offering the best balance between speed and stability in convergence.

\section{Conclusion}
The integration of Enhanced-GD, automatic-, and numerical differentiation techniques with computational graphs significantly advances power flow analysis. This method has proven its efficacy through rigorous simulations, adeptly handling large-scale ACPF challenges, especially in complex systems. The adaptive learning rate of Enhanced-GD effectively minimizes errors and expedites convergence, significantly outperforming traditional approaches. Additionally, computational graphs not only optimize calculations but also improve the visualization and interpretation of intricate, non-convex power systems. This cutting-edge strategy equips analysts with a robust toolkit for non-convex power flow analysis, ensuring resilience and precision amidst evolving grid dynamics.

\end{document}